\documentclass[12pt,a4paper]{article}
\usepackage{latexsym}
\usepackage{amssymb}
\newtheorem{lem}{Lemma}
\newtheorem{prop}{Proposition}
\newtheorem{thm}{Theorem}
\newenvironment{pf}{\textbf{Proof\ }}{\hfill$\Box$\smallskip}
\newcommand{\la}{\lambda}
\newcommand{\om}{\omega}
\newcommand{\ca}{\mathcal{A}}
\newcommand{\co}{\mathcal{O}}
\newcommand{\ct}{\mathcal{T}}
\newcommand{\cv}{\mathcal{V}}
\newcommand{\cx}{\mathcal{X}}
\newcommand{\cd}{\bullet}
\newcommand{\Na}{\uparrow}
\newcommand{\Sa}{\downarrow}
\newcommand{\Ea}{\rightarrow}
\newcommand{\Wa}{\leftarrow}
\title{Alternating sign matrices and tournaments}
\author{Robin Chapman\\
School of Mathematical Sciences\\ University of Exeter\\
Exeter, EX4 4QE, UK\\ \texttt{rjc@maths.ex.ac.uk}}
\date{18 October 2000}
\begin{document}
\maketitle

\section{Introduction}

An alternating sign matrix of order $n$ is an $n$-by-$n$ matrix
with entries in $\{-1,0,1\}$ such that each row and column has sum 1 and
the nonzero entries in each row and column alternate in sign. For example
\begin{equation}\label{exASM}
\left(\begin{array}{rrrrr}
0& 0&0& 1&0\\
0& 0&1& 0&0\\
0& 1&0&-1&1\\
1&-1&0& 1&0\\
0& 1&0& 0&0
\end{array}\right).
\end{equation}
The problem of enumerating alternating sign matrices has led to many
deep results and difficult conjectures. For a history
of alternating sign matrices see~\cite{BrPC}.

There are several combinatorial and physical structures equivalent to
alternating sign matrices. One of the most fertile has been the
square-ice model originally used in statistical mechanics. The square
ice graph corresponding to (\ref{exASM}) is
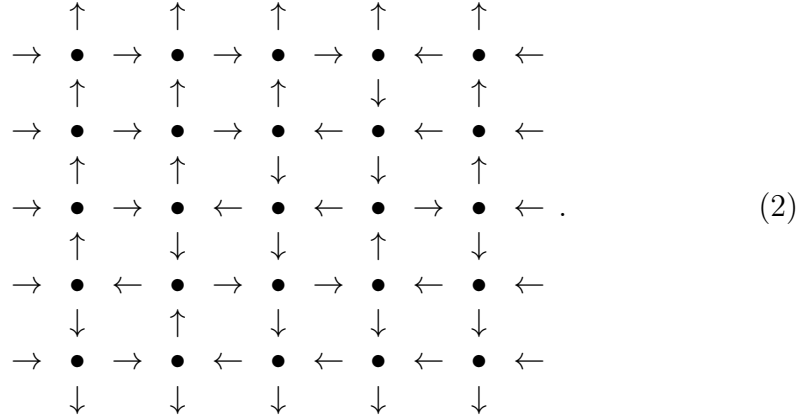
\begin{equation}
\begin{array}{ccccccccccc}
   &\Na&   &\Na&   &\Na&   &\Na&   &\Na&   \\
\Ea&\cd&\Ea&\cd&\Ea&\cd&\Ea&\cd&\Wa&\cd&\Wa\\
   &\Na&   &\Na&   &\Na&   &\Sa&   &\Na&   \\
\Ea&\cd&\Ea&\cd&\Ea&\cd&\Wa&\cd&\Wa&\cd&\Wa\\
   &\Na&   &\Na&   &\Sa&   &\Sa&   &\Na&   \\
\Ea&\cd&\Ea&\cd&\Wa&\cd&\Wa&\cd&\Ea&\cd&\Wa\\
   &\Na&   &\Sa&   &\Sa&   &\Na&   &\Sa&   \\
\Ea&\cd&\Wa&\cd&\Ea&\cd&\Ea&\cd&\Wa&\cd&\Wa\\
   &\Sa&   &\Na&   &\Sa&   &\Sa&   &\Sa&   \\
\Ea&\cd&\Ea&\cd&\Wa&\cd&\Wa&\cd&\Wa&\cd&\Wa\\
   &\Sa&   &\Sa&   &\Sa&   &\Sa&   &\Sa&
\end{array}.
\end{equation}
This directed graph has inward edges at the left and right and
outward edges at the top and bottom. Both the in-degree and out-degree
of all its degree four vertices equal two. These vertices are of
six types, which we label according to their incident in-edges.
\emph{Horizontal} vertices have two horizontal in-edges,
\emph{vertical} vertices have two vertical in-edges and
\emph{northwest} vertices have in-edges pointing in the north
(up) and west (left) directions. We define \emph{southwest},
\emph{northeast} and \emph{southeast} vertices similarly. The
correspondence between alternating sign matrices and square-ice
graphs is characterized as follows: horizontal vertices correspond
to entries~1, vertical vertices to entries $-1$ and all
other vertices to entries~0.

Bressoud \cite{Br3ASM} has asked for a bijective proof of the following
identity
\begin{equation}\label{BrId}
\sum_{T\in\ct_n}\la^{U(T)}\prod_{i=1}^nx_i^{\om(i)}
=\sum_{A\in\ca_n}\la^{SW(A)}(1+\la)^{V(A)}
\prod_{i=1}^nx_i^{SW_i(A)+SE_i(A)+V_i(A)}
\end{equation}
due essentially to Robbins and Rumsey~\cite{RR}.
Here $\ca_n$ is the set of $n$-by-$n$ alternating sign matrices,
$SW(A)$, $SE(A)$ and $V(A)$ denote respectively the numbers of
southwest, southeast and vertical vertices in the square-ice graph
corresponding to $A\in\ca_n$, and $SW_i(A)$, $SE_i(A)$ and $V_i(A)$
denote the numbers of vertices of the appropriate type in the
$i$-th column of this graph. Also $\ct_n$ is the set of tournaments
on the set $\{1,2,\ldots,n\}$, that is the set of orientations on
the complete graph with vertex set $\{1,2,\ldots,n\}$.
Given $T\in\ct_n$ we denote by $\om(i)$ the out-degree of the vertex $i$
in~$T$, and we let $U(T)$ be the number of upsets in~$T$, that is the number
of edges oriented from $i$ to $j$ with $i>j$.

We give a bijective proof of the following refinement of (\ref{BrId}):
\begin{equation}\label{RefId}
\prod_{1\le i<j\le n}(x_i+y_j)
=\sum_{A\in\ca_n}\prod_{i=1}^n x_i^{SE_i(A)}y_i^{SW_i(A)}(x_i+y_i)^{V_i(A)}.
\end{equation}
This equation is equivalent to \cite[Theorem~3.1]{C}.
Putting $y_i=\la x_i$ in (\ref{RefId}) gives (\ref{BrId}) since the factor
$x_i+\la x_j$ on the left corresponds to the edge between $i$ and~$j$. This
edge is either oriented from $i$ to $j$ and so contributes $x_i$ to the
product, or is oriented from $j$ to $i$ and so contributes $\la x_j$ to the
product. The latter case corresponds to an upset and so the correct
power of $\la$ appears in the product.

\section{Triangles}

We use an alternative combinatorial model for alternating sign matrices,
namely complete monotone triangles. A complete monotone triangle
of order $n$ is
a triangular array of integers with $n$ rows and $k$ entries in the
$k$-th row. Each row is strictly increasing
and the final row contains the entries $1,2,\ldots,n$ in
ascending order. Each entry not in the final row lies weakly
between the entries immediately to the left and right in the row below.
Here is an example of a complete monotone triangle of order 5:
\begin{equation}\label{exCMT}
\begin{array}{ccccccccc}
 & & & &4& & & & \\
 & & &3& &4& & & \\
 & &2& &3& &5& & \\
 &1& &3& &4& &5& \\
1& &2& &3& &4& &5
\end{array}.
\end{equation}

To produce a complete monotone triangle from an alternating sign matrix
we form a new matrix by taking its cumulative column sums. For instance
(\ref{exASM}) has cumulative column sums
\begin{equation}\label{colsum}
\left(\begin{array}{ccccc}
0&0&0&1&0\\
0&0&1&1&0\\
0&1&1&0&1\\
1&0&1&1&1\\
1&1&1&1&1
\end{array}\right).
\end{equation}
This new matrix has $k$ 1s in the $k$-th row with the remaining entries~0.
We form a triangle by writing in the $k$-th row of the triangle,
the positions of the $k$ 1s in the column sum matrix. Applying this to
(\ref{colsum}) gives (\ref{exCMT}). This gives a bijection between
alternating sign matrices and complete monotone triangles of the same
order.

The entries not in the bottom row of a complete monotone triangle fall into
three types. Consider an entry $j$ and let $i$ and $k$ be its neighbours
in the row below, so that
$$\begin{array}{ccc}
 &j& \\
i& &k
\end{array}$$
is part of the triangle. Then either
\begin{description}
\item{(i)} $i=j<k$,
\item{(ii)} $i<j=k$, or
\item{(iii)} $i<j<k$.
\end{description}
Entries of types (i), (ii) and (iii) correspond respectively
to southeast, southwest and vertical vertices in column $j$ of the
associated square-ice graph. We introduce a refinement of the
notion of complete monotone triangle, namely that of \emph{oriented complete
monotone triangle}. An oriented complete monotone triangle is a
complete monotone triangle with an orientation (left or right)
given to each entry not in the bottom row. Entries of type (i)
must be oriented to the right and entries of type (ii) must be oriented
to the left. Entries of type (iii) may be oriented in either direction.
We shall depict an oriented complete monotone triangle by replacing
each right-oriented entry $j$ by a symbol $x_j$ and
each left-oriented entry $j$ by a symbol $y_j$. Thus in the
configuration
$$\begin{array}{cc}
x_j&   \\
   &z_k
\end{array}$$
where $z_k$ denotes either $x_k$, $y_k$ or $k$ (in the bottom row)
we must have $j<k$ while in
$$\begin{array}{cc}
   &y_j\\
z_i&
\end{array}$$
we must have $i<j$.

As an example, one of the oriented
complete monotone triangles arising from (\ref{exCMT}) is
\begin{equation}\label{exOCMT}
\begin{array}{ccccccccc}
 &   &   &   &y_4&   &  &    &\\
 &   &   &y_3&   &x_4&  &    &\\
 &   &y_2&   &x_3&   &y_5&   &\\
 &x_1&   &y_3&   &y_4&   &y_5&\\
1&   &  2&   &  3&   &  4&   &5
\end{array}.
\end{equation}

As the type (iii) vertices in a complete monotone triangle correspond
to the $-1$ entries in the alternating sign matrix, an alternating sign
matrix with $r$ $-1$s gives rise to $2^r$ oriented complete monotone
triangles. Denote the set of oriented complete monotone triangles
of order $n$ by~$\co_n$.

We define the \emph{weight} of an oriented complete monotone triangle
to be the product of the entries in the rows other than the last. Thus
the weight of the triangle in (\ref{exOCMT}) is
$x_1 y_2 x_3 y_3^2 x_4 y_4^2 y_5^2$.
It is clear that the right side of (\ref{RefId}) is the sum of the weights of
all elements of $\co_n$.

We define the weight of a tournament $T\in\ct_n$ as the products
of the weights of its edges where the weight of an edge oriented from $i$ to
$j$ is $x_i$ if $i<j$ and $y_i$ if $i>j$. Then the left side of (\ref{RefId})
is the sum of the weights of all elements of $\ct_n$.
We shall prove the following theorem.

\begin{thm}\label{main}
There is a weight-preserving bijection between $\co_n$ and $\ct_n$.
\end{thm}

The equation (\ref{RefId}) is an immediate corollary. The remainder of the
paper is devoted to the construction of the weight-preserving bijection
$\Phi:\co_n\to\ct_n$ and its inverse $\Psi:\ct_n\to\co_n$.

We first indicate how to display a tournament in a triangle in a similar
fashion to an oriented complete monotone triangle.
We introduce new indeterminates
$a_{ij}$ and $b_{ij}$ where $1\le i<j\le n$. For $1\le k<n$, let
$\cv_k=\{a_{i,j},b_{i,j}:j-i=k\}$. Also let
$\cx=\{x_i,y_j:1\le i,j\le n\}$. Each tournament corresponds to
a set of $\frac12n(n-1)$ of indeterminates in $\bigcup_{k=1}^{n-1}\cv_k$
where we take $a_{ij}$
when there is an edge oriented from $i$ to $j$ and we take $b_{ij}$
when there is an edge oriented from $j$ to~$i$. We arrange these
symbols in $n-1$ rows of a triangle where we put those
entries in $\cv_{n-j}$ in the $j$-th row in `ascending order'.
The $i$-th entry in row $j$ thus has subscript ${}_{i\,(i+n-j)}$.
We fill the $n$-th row with the numbers $1,2,\ldots,n$
in ascending order. For example,
\begin{equation}\label{exTour}
\begin{array}{ccccccccc}
 &      &      &      &a_{15}&      &      &      &\\
 &      &      &b_{14}&      &b_{25}&      &      &\\
 &      &b_{13}&      &b_{24}&      &b_{35}&      &\\
 &b_{12}&      &b_{23}&      &a_{34}&      &a_{45}&\\
1&      &     2&      &     3&      &     4&      &5
\end{array}.
\end{equation}

We now identify tournaments with triangles of this form.
We define the bijections $\Phi$ and $\Psi$ by taking a triangle of
one type and converting it by a series of moves to one of the other type.
We pass through a series of intermediate triangles which in general
contain variables $x_i$ and $y_i$ and variables $a_{ij}$ and $b_{ij}$.
We define a \emph{triangle} of order $n$ as an $n$-rowed triangle
with $k$ entries in the $k$-th row and
with $1,2,\ldots,n$ in the last row and entries drawn from
$\cx\cup\bigcup_{k=1}^{n-1}\cv_k$ in the first $(n-1)$-rows.
Each of our primitive moves either permutes the entries of
the triangle or replaces an $a_{ij}$ by $x_i$ or an $x_i$ by an $a_{ij}$
or a $b_{ij}$ by $y_j$ or a $y_j$ by a $b_{ij}$. We define
the \emph{weight} of a triangle as the product of the weights of its entries
above its bottom row
where $x_i$, $y_j$, $a_{ij}$ and $b_{ij}$ have weights $x_i$, $y_j$,
$x_i$ and $y_j$ respectively. This weighting agrees with that of oriented
monotone triangles and tournaments, and will be
preserved under each primitive move.

We interpret $x_j$ as an outgoing edge from vertex $j$ to
some as yet unknown larger vertex. Similarly $y_j$ represents
an outgoing edge from vertex $j$ to an unknown smaller vertex. Also
$a_{ij}$ denotes an edge oriented from $i$ to $j$ ($i<j$)
and $b_{ij}$ denotes an edge oriented from $j$ to $i$ ($i<j$).
At intermediate stages of our algorithm, the $a_{ij}$s and $b_{ij}$s
that occur are still provisional, as for instance $a_{ij}$ may be replaced
by $x_i$ and then by $a_{ik}$ for $k\ne j$.

For convenience we suppose that our triangles have an empty ``row 0'' above
the first row. We also use the symbols $z_i$ to denote indifferently $x_i$,
$y_i$ or $i$ (appearing in the final row) and $c_{ij}$ to denote indifferently
$a_{ij}$ or $b_{ij}$.

We need to define the notion of an \emph{admissible} triangle. A
triangle with no entry in $\bigcup_{k=1}^{n-1}\cv_k$ will be
admissible if and only if it lies in $\co_n$ and a
triangle with no entry in $\cx$ will be
admissible if and only if it lies in $\ct_n$. At each stage of our
algorithm, admissibility will be preserved.

A triangle is admissible if and only if the following conditions
hold:
\begin{description}
\item{(i)} it has an admissible ranking,
\item{(ii)} it has the monotone diagonal property,
\item{(iii)} it has the monotone row property.
\end{description}

We explain the conditions in turn. We say that a row of a triangle
is an $\cx$-row if all its entries lie in $\cx$ and that it is a
is a $\cv$-row if all its entries lie in some $\cv_k$.
(We regard the notional row 0 as both a $\cx$-row and a $\cv$-row).
An \emph{admissible} ranking of a triangle is an assignment of non-negative
integers, called ranks, to each row of the triangle (including row 0)
with the following properties:
\begin{itemize}
\item the $n$-th row has rank 0 and the $(n-1)$-th row has rank 1,
\item if the $r$-th row ($r>0$) has rank $k$ then either row $r$ is
a $\cx$-row and row $r-1$ also has rank $k$, or row $r$ has all its entries
in $\cv_k$ and row $r-1$ has rank $k+1$.
\end{itemize}
It follows that the ranks in an admissible ranking decrease weakly
as one descends the triangle and that the ranks of the $\cv$-rows
form the sequence $d,d-1,\ldots,2,1$ for some $d\ge0$. The rank of
a $\cv$-row with entries in $\cv_k$ is~$k$.

If a triangle has an admissible ranking then it is plain that
the ranking is uniquely determined. We call such a triangle
a \emph{ranked} triangle.

To explain the monotone diagonal and monotone row properties we need
the notion of the left and right values of an entry in a ranked triangle.
Let $u$ be an entry in a ranked triangle. It is assigned a \emph{left value}
$l(u)$ and a \emph{right value} $r(u)$. Note that these values will
in general depend on the rank of the row in which $u$ lies, not simply
on the indeterminate $u$ itself.
When $u=a_{ij}$ or $b_{ij}$ we define $l(u)=i$ and $r(u)=j$.
If $u=x_i$ lies in a row with rank $t$ we define $l(u)=i$
and $r(u)=i+t$ while if $u=y_j$ lies in a row with rank $t$ we define
$l(u)=j-t$ and $r(u)=j$. Finally if $u=i$ is a number in the last row we let
$l(u)=l(v)=i$.
In all cases $r(u)-l(u)$ is the rank of the row containing~$u$.

To define the monotone diagonal property consider three entries
in a ranked triangle appearing in the following arrangement
\begin{equation}\label{uvw}
\begin{array}{ccc}
 &v& \\
u& &w
\end{array}.
\end{equation}
An arrangement (\ref{uvw}) has the \emph{monotone diagonal property} if
\begin{itemize}
\item $l(u)\le l(v)$ and $r(v)\le r(w)$,
\item if $u=y_j$ for some $j$ then $l(u)<l(v)$, and
\item if $w=x_j$ for some $j$ then $r(v)<r(w)$.
\end{itemize}
A ranked triangle has the monotone diagonal property
if all arrangements (\ref{uvw}) within have the monotone diagonal
property.
The first of these conditions states that the left values increase
weakly when traversing northeast diagonals and the right values increase
weakly when traversing southeast diagonals.
If the second condition is relevant, then the two rows represented
in (\ref{uvw}) have the same rank, and we can
rewrite the condition as $j<r(v)$. Similarly the final condition
is equivalent to $l(v)<j$.

If $T$ is a ranked triangle having the monotone diagonal
property then the
the left value of the $j$-th entry in each row is
at least $j$ and that the right value of the $k$-th entry from the right
in each row is at most $n+1-k$.

We say that an arrangement (\ref{uvw})
has the \emph{monotone row property} if
\begin{itemize}
\item $l(u)<l(w)$.
\end{itemize}
We say that a ranked triangle
has the \emph{monotone row property} if
every arrangement (\ref{uvw}) within has the monotone row property.

\begin{lem}
A triangle of order $n$ all of whose rows but the last are $\cx$-rows
is admissible if and only if it lies in $\co_n$.
\end{lem}
\begin{pf}
Suppose that $T\in\co_n$. Then $T$ clearly has an admissible ranking
with each row but the last having rank~1. Let us consider the
arrangement~(\ref{uvw}). We divide into cases according to whether
$u$, $v$ and $w$ are $x$s or $y$s or they lie in the last row.
The cases where $v=y_j$ are
similar to those where $v=x_j$ so we look at these in detail:
\begin{description}
\item{(i)} $u=x_i$, $w=x_k$. Then $i\le j<k$.
Hence $l(u)=i\le j=l(v)$ and $r(v)=j+1<k+1=r(w)$.
Also $l(u)=i<k=l(w)$.

\item{(ii)} $u=x_i$, $w=y_k$. Then $i\le j<k$.
Hence $l(u)=i\le j=l(v)$ and $r(v)=j+1\le k=r(w)$.
We consider the entry between $u$ and $w$ in the next row:
$$\begin{array}{ccc}
   &x_j&   \\
x_i&   &y_k\\
   &z_l&
\end{array}.$$
We have $i<l<k$ so that $i+1<k$ and so $l(u)=i<k-1=l(w)$.

\item{(iii)} $u=y_i$, $w=x_k$. Then $i\le j<k$.
Hence $l(u)=i-1<j=l(v)$ and $r(v)=j<k=r(w)$. Also
$l(u)+1=i<k=l(w)$.

\item{(iv)} $u=y_i$, $w=y_k$. Then $i\le j<k$. Hence $l(u)=i-1<j=l(v)$
and $l(v)=j\le k-1=l(w)$. Also $l(u)=i-1<k-1=l(w).$

\item{(v)} $u=i$, $w=k$ (in the last row). Then $i\le j<k$. Hence
$l(u)=i\le j=l(v)$ and $r(v)=j+1\le k=r(w)$.

\end{description}
In each case, the monotone diagonal and monotone row
properties are satisfied. Hence $T$ is an admissible triangle.

Now suppose that $T$ is an admissible triangle with each row but the
last an $\cx$-row. Then all rows but the last are ranked 1.
Again consider the arrangement (\ref{uvw}) and divide into cases as above.
As before we consider in detail only the cases where $v=x_j$
as the cases where $v=y_j$ are similar.

\begin{description}
\item{(i)} $u=x_i$, $w=x_k$. Then $i=l(u)\le l(v)=j$
and $j=r(v)-1<r(w)-1=k$. Hence $i<k$.

\item{(ii)} $u=x_i$, $w=y_k$. Then $i=l(u)\le l(v)=j$
and $j=r(v)-1\le r(w)-1<r(w)=k$. Hence $i<k$.

\item{(iii)} $u=y_i$, $w=x_k$. Then $i=l(u)+1<l(v)+1=j+1$
and so $i\le j$. Also $j=r(v)-1<r(w)-1=k$. Hence $i<k$.

\item{(iv)} $u=y_i$, $w=y_k$. Then $i=l(u)+1<l(v+1)=j+1$
and so $i\le j$. Also $j=r(v)-1\le r(w)-1=k-1$ and so $j<k$. Hence $i<k$.

\item{(v)} $u=i$, $w=k$ (in the last row). Then $i=l(u)\le l(v)=j$
and $j=r(v)-1\le r(w)-1=k-1$ and so $j<k$.
\end{description}
It follows that $T$ is an oriented complete monotone triangle.
\end{pf}

\begin{lem}
A triangle of order $n$ all of whose rows but the last are $\cv$-rows
is admissible if and only if it lies in $\ct_n$.
\end{lem}
\begin{pf}
If the triangle $T$ with no $\cx$-rows is admissible then
the rank of the $k$-th row must be $t=n-k$. Its entries
must contain each of the $k$ possible subscripts of the form ${}_{i\,i+t}$
in increasing order, by the monotone row property. It must therefore be
a tournament.

Conversely suppose that $T$ is a tournament. In any arrangement (\ref{uvw}) if
$l(u)=i$ then $r(u)=i+t$ where $t$ is the rank of the lower row and
then $l(v)=i$, $r(v)=i+t+1$, $l(w)=i+1$ and $r(w)=i+t+1$. It is then
plain that the monotone diagonal and monotone row properties hold, so
that $T$ is admissible.
\end{pf}

\section{The algorithm}

We define two operations on admissible triangles. Each operates on a pair
of consecutive rows: rows $d-1$ and $d$ where $1\le d<n$. The operation
of \emph{raising} can be applied when row $d$ is a $\cx$-row and
and row $d-1$ is a $\cv$-row (including the case when $d-1=0$).
Rows $d-1$ and $d$ then have the same
rank~$t$. We start by making the following
swaps in rows $d-1$ and $d$ whenever possible
$$\begin{array}{cc}
   &b_{jk}\\
x_j&
\end{array}\longrightarrow
\begin{array}{cc}
      &x_j\\
b_{jk}&
\end{array}
$$
and
$$\begin{array}{cc}
a_{jk}&   \\
      &y_k
\end{array}\longrightarrow
\begin{array}{cc}
y_k&      \\
   &a_{jk}
\end{array}.
$$
We then replace any $a_{jk}$
remaining in row $d-1$ by $x_j$ and any $b_{jk}$
remaining in row $d-1$ by $y_k$. We then replace all $x_i$ in row $d$
by $a_{i\,i+t}$ and all $y_j$ in row $r$ by $b_{j-t\,j}$. For instance
$$\begin{array}{ccccccc}
   &b_{13}&   &a_{57}&   &a_{68}&   \\
x_1&      &y_5&      &y_7&      &x_7
\end{array}$$
first becomes
$$\begin{array}{ccccccc}
      &x_1&   &y_7&      &a_{68}&   \\
b_{13}&   &y_5&   &a_{57}&      &x_7
\end{array}$$
and finally becomes
$$\begin{array}{ccccccc}
      &x_1&      &y_7&      &x_6&        \\
b_{13}&   &b_{35}&   &a_{57}&   &a_{79}.
\end{array}$$
It is not immediately obvious that when we introduce an
$a_{i\,i+t}$ or a $b_{j-t\,j}$ that $i+t\le n$ or $j-t\ge1$.
For the moment we shall ignore this difficulty. Later
we will show that raising preserves admissibility. As the
left and right values of each entry in an admissible triangle
lie between 1 and $n$ inclusive, then we must have
$i+t\le n$ or $j-t\ge1$ as appropriate.

We denote this raising operation by $R_d$. The rank of row $d-1$ has now
been raised to $t+1$ while the ranks of the remaining rows are left
unchanged. Also the weight of the new triangle is the same as that of
the old.

\begin{prop}
When applied to an admissible triangle, the raising operator $R_d$
produces an admissible triangle.
\end{prop}
\begin{pf}
Let $T$ be an admissible triangle with row $d-1$ and $d$ having
rank $t$ and with row $d-1$ a $\cv$-row and row $d$ an $\cx$-row.
Let $T'=R_d(T)$.

First note that the operation of raising does not alter the left and
right values of the entries in row~$d$. In row $d-1$ suppose the
entry $v$ is replaced by~$v'$. If
$v'=x_i$ then $l(v')=l(v)$ and $r(v')=r(v)+1$, while if
$v'=y_j$ then $r(v')=r(v)-1$ and $r(u')=r(u)$.

Suppose an arrangement (\ref{uvw}) in $T$ becomes
\begin{equation}\label{uvwdash}
\begin{array}{ccc}
  &v'&  \\
u'&  &w'
\end{array}
\end{equation}
in~$T'$.
We need to check that (\ref{uvwdash}) has the
monotone diagonal and monotone row properties. Let the lower
row in (\ref{uvwdash}) be row $e$. The entries
in (\ref{uvw}) and (\ref{uvwdash}) will be the same unless
$e=d-1$, $d$ or $d+1$. Also we only need to show the monotone row
property when $e=d-1$ or~$d$.

We can swiftly dispose of the monotone diagonal property when $e=d+1$.
In this case $u=u'$ and $w=w'$. Also $l(v)=l(v')$ and $r(v)=r(v')$.
The monotone diagonal property for (\ref{uvwdash}) follows immediately
from that for (\ref{uvw}).

Now consider the case where $e=d-1$. In this case $l(u')=l(u)$
or $l(u)-1$ and $r(w')=r(w)$ or $r(w)+1$. Also $v=v'$.
Thus $l(u')\le l(u)\le l(v)=l(v')$ and $r(v')=r(v)\le r(w)\le r(w')$.
If $u'=y_j$, then $l(u')=l(u)-1<l(v)=l(v')$. Similarly if
$w'=x_j$, then $r(v')=r(v)<r(w)+1=r(w')$. Also $l(u')<l(w')$
unless $l(w)=l(u)+1$, $u'=x_i$ and $w'=y_j$ for some $i$ and~$j$.
Then $l(u)=i$, $r(u)=j-1$, $l(w)=i+1$, $r(w)=j$ and $j-i=t+1$.
Consider the entry in $T$ below $u$ and~$w$:
$$\begin{array}{ccc}
 &v& \\
u& &w\\
 &q&
\end{array}.$$
From the monotone diagonal property of~$T$, $q$ must be either
$x_{i+1}$ or $y_{j-1}$. Suppose the former. Then $w$ cannot be
$b_{i+1\,j}$ since that would be swapped with $q$ yielding $w'=x_{i+1}$.
Thus $w=a_{i+1,j}$ and for $w'$ to be $y_j$ we need the entry following
$q$ in row $d$ in $T$ to be $y_j$:
$$\begin{array}{cccc}
u&       &a_{i+1\,j}&   \\
 &x_{i+1}&          &y_j
\end{array}.$$
But then in row~$d$, $l(x+1)=i+1$ and $r(x+1)=j-t=i+1$, violating the
monotone row property for~(\ref{uvw}). A similar
argument yields a contradiction when $q=y_{j-1}$. Hence (\ref{uvwdash})
also has the monotone diagonal and monotone row properties.

Finally suppose that $d=e$. Now $l(u)=l(u')$ and $l(v')\in\{l(v)-1,l(v)\}$.
Thus $l(u')\le l(v')$ unless $l(u)=l(v)$ and $v'=y_j$. Then $u=x_{j-t}$
or $y_j$ and $v=a_{j-t\,j}$ or $v=b_{j-t\,j}$. If $u=y_j$ then
$l(u)=j-t=l(v)$ contrary to the monotone diagonal property for~(\ref{uvw}).
Thus $u=x_{j-t}$. If $v=b_{j-t\,j}$ then $u$ and $v$ are swapped giving
$v'=x_{j-t}$. Thus $v=a_{j-t\,j}$. As $v'=y_j$ then $y_j$ must succeed
$x_{j-t}$ in row $d$ of~$T$:
$$\begin{array}{ccc}
       &a_{j-t\,j}&   \\
x_{j-t}&          &y_j
\end{array}.$$
But then in row $d$ of~$T$, $l(x_{j-t})=j-t=l(y_j)$ contrary to the
monotone row property of $T$. Thus $l(u')\le l(v')$.
Similarly $r(v')\le r(w')$. Thus (\ref{uvw}) has the monotone diagonal
property. As row $d$ in $T'$ has the same left and right values as the
row $d$ of~$T$, then (\ref{uvwdash}) has the monotone
row property. We conclude that the triangle $T'$ is admissible.
\end{pf}

The inverse of raising is \emph{lowering}. This time suppose that row $d$
($1\le d<n$) is a $\cv$-row of rank $t$ and that
row $d-1$ (including the case $d-1=0$) is a $\cx$-row.
Then row $d-1$ has rank $t+1$. We make the swaps
$$\begin{array}{cc}
      &x_j\\
b_{jk}&
\end{array}\longrightarrow
\begin{array}{cc}
   &b_{jk}\\
x_j&
\end{array}
$$
and
$$\begin{array}{cc}
y_k&      \\
   &a_{jk}
\end{array}\longrightarrow
\begin{array}{cc}
a_{jk}&   \\
      &y_k
\end{array}
$$
where possible in rows $d-1$ and $d$. We then replace the
$a_{jk}$ and $b_{jk}$ in row $d-1$ by $x_j$ and $y_k$ and the
$x_j$ and $y_k$ in row $d$ by
$a_{j\,j+t}$ and $b_{k-t\,k}$. Again, it will follow from the fact that
lowering preserves admissibility that $j+t\le n$ or $k-t\ge1$ as appropriate.
The new triangle has an admissible ranking.
All the ranks are unchanged save for row $d-1$ whose rank is lowered from
$t+1$ to~$t$. We denote this lowering operator by $L_d$.

\begin{prop}
When applied to an admissible triangle, the lowering operator $L_d$
produces an admissible triangle.
\end{prop}
\begin{pf}
Let $T$ be an admissible triangle with row $d-1$ an $\cx$-row and row
$d$ a $\cv$-row of rank~$t$. Then row $d-1$ of $T$ has rank $t+1$.
Let $T'=L_d(T)$.

The lowering operation does not alter the left and right values in row~$d$.
Suppose that the entry $v$ in row $d-1$ is replaced by~$v'$.
If $v=x_j$ then $l(v')=l(v)$ and $r(v')=r(v)-1$ while if $v=y_j$ then
$l(v')=l(v)+1$ and $r(v')=r(v)$.

Suppose we have an arrangement (\ref{uvw}) having lower row $e$, becoming
(\ref{uvwdash}) after lowering. They are identical unless $e\in\{d-1,d,d+1\}$.

When $e=d+1$ the monotone row property of (\ref{uvwdash}) is immediate. Also
the monotone diagonal property follows as the values in row $d$ and
the entries in row $d+1$ are unchanged.

Suppose that $e=d-1$. Row $d-1$ in $T'$ is a $\cv$-row. We have
$l(u')\in\{l(u),l(u)+1\}$ and $l(v')=l(v)$ so that $l(u')\le l(v')$
unless $l(u)=l(v)$ and $l(u')=l(u)+1$. In this case $u=y_j$ and
so $l(u)<l(v)$ contrary to hypothesis. Hence $l(u')\le l(v')$ and similarly
$r(v')\le r(w')$. Thus (\ref{uvwdash})
has the monotone diagonal property.
As also $l(w')\in\{l(w),l(w)+1\}$ and $l(u)<l(w)$ then $l(u')<l(v')$
unless $l(u)+1=l(w)$, $l(u')=l(u)+1$ and $l(w')=l(w)$. In this case
$u=y_i$ and $w=x_j$ so that by the monotone diagonal property for (\ref{uvw})
we have $l(u)<l(v)$ and $r(v)<r(w)$. Consequently $l(v)<l(w)$ and
so $l(w)>l(u)+1$ contrary to hypothesis. Hence (\ref{uvwdash}) also has
the monotone row property.

Finally suppose that $e=d$. As $l(v')\in\{l(v),l(v)+1\}$ then
$l(u')=l(u)\le l(v)\le l(v')$. Similarly $r(v')\le r(v)\le r(w)=r(w')$.
Suppose that $u'=y_i$. We must show that $l(u')<l(v')$. This is valid unless
$l(u')=l(v')$. If so their common value is $i-t$. As $l(v')=i-t$
and $r(v')=t$ then $v=x_{i-t}$ or $y_i$. Also $u=a_{i-t\,i}$
or $b_{i-t\,i}$. We cannot have $v=y_i$ for then $l(v)=i-t-1$
and $l(u)=i-t$ so that $l(u)\not\le l(v)$. But if
$v=x_{i-t}$ and $u=b_{i-t\,i}$ then $u$ and $v$ are swapped so that
$u'=x_{i-t}$. If $v=x_{i-t}$ and $u=a_{i-t\,i}$ the only way that $u'$
can be equal to $y_i$ is if the entry preceding $u$ in the $(d-1)$-th
row of $T$ is $y_i$ so that we get
$$\begin{array}{cccc}
y_i&          &x_{i-t}& \\
   &a_{i-t\,i}&       &w
\end{array}.$$
But in this case the left values of the entries in row $d-1$ are
$i-t-1$ and $i-t$ which violates the the condition that $T$
has the monotone row property. Hence $l(u')<l(v')$.
Similarly if $v'=x_j$ then $r(v')<r(w')$ and so (\ref{uvwdash}) has the
monotone diagonal property.

As the values of the entries in row $d$ do not change, then (\ref{uvwdash})
has the monotone row property.
We conclude that the triangle $T'$ is admissible.
\end{pf}

\begin{lem}
The operations $R_d$ and $L_d$ are inverse: if we operate on an
admissible triangle by $R_d$ and then $L_d$ we regain the initial
triangle. Similarly if we follow $L_d$ by $R_d$.
\end{lem}
\begin{pf}
Suppose that $T$ is an admissible triangle and that $T'=R_d(T)$ is defined.
When we raise $T$ by $R_d$ we first swap some pairs of entries and
then when we lower $T'$ by $L_d$ these entries are swapped
back. If we can show that no
other pairs of entries are swapped when lowering we are done, since
it is plain that the remaining entries return to their original
state. Suppose after raising we have a configuration
$$\begin{array}{cc}
      &x_j\\
b_{jk}&
\end{array}$$
which has not arisen from swapping the entries. This must have come
from the configuration
$$\begin{array}{cc}
   &a_{jk}\\
y_k&
\end{array}$$
which is banned by the monotone diagonal property. Similarly a configuration
$$\begin{array}{cc}
y_k&      \\
   &a_{jk}
\end{array}$$
cannot arise except by swapping. Hence raising followed by lowering gives the
identity. Similarly lowering followed by raising gives the identity.
\end{pf}

We can now describe the bijections $\Phi$ and~$\Psi$. We define
$$\Phi=(R_1\circ R_2\circ R_3\circ\cdots\circ R_{n-1})\circ\cdots
\circ(R_1\circ R_2\circ R_3)\circ(R_1\circ R_2)\circ(R_1).$$
It is easy to see that $\Phi(T)$ is well-defined for $T\in\co_n$
and that $\Phi(T)\in\ct_n$. Its inverse $\Psi$ is then
$$(L_1)\circ(L_2\circ L_1)\circ(L_3\circ L_2\circ L_1)\circ\cdots\circ
(L_{n-1}\circ\cdots\circ L_3\circ L_2\circ L_1).$$
As $\Phi$ and $\Psi$ are composites of weight-preserving maps, they too
are weight-preserving.
This completes the proof of Theorem~\ref{main}

\section{An example}

As an example let us apply $\Phi$ to the oriented complete monotone
triangle~(\ref{exOCMT}). We indicate the rank of each row in brackets at the
end, including the rank of the empty row 0.
Here is (\ref{exOCMT}) in this notation:
$$\begin{array}{cccccccccc}
 &   &   &   &   &   &   &   & &\qquad[1]\\
 &   &   &   &y_4&   &   &   & &\qquad[1]\\
 &   &   &y_3&   &x_4&   &   & &\qquad[1]\\
 &   &y_2&   &x_3&   &y_5&   & &\qquad[1]\\
 &x_1&   &y_3&   &y_4&   &y_5& &\qquad[1]\\
1&   &  2&   &  3&   &  4&   &5&\qquad[0]
\end{array}.$$
We first apply the operator $R_1$ giving
$$\begin{array}{cccccccccc}
 &   &   &   &      &   &   &   & &\qquad[2]\\
 &   &   &   &b_{34}&   &   &   & &\qquad[1]\\
 &   &   &y_3&      &x_4&   &   & &\qquad[1]\\
 &   &y_2&   &   x_3&   &y_5&   & &\qquad[1]\\
 &x_1&   &y_3&      &y_4&   &y_5& &\qquad[1]\\
1&   &  2&   &     3&   &  4&   &5&\qquad[0]
\end{array}.$$
Next apply $R_2$:
$$\begin{array}{cccccccccc}
 &   &   &      &     &      &   &  & &\qquad[2]\\
 &   &   &      &  y_4&      &   &  & &\qquad[2]\\
 &   &   &b_{23}&     &a_{45}&   &  & &\qquad[1]\\
 &   &y_2&      &  x_3&      &y_5&  & &\qquad[1]\\
 &x_1&   &   y_3&     &   y_4&  &y_5& &\qquad[1]\\
1&   &  2&      &    3&      &  4&  &5&\qquad[0]
\end{array}.$$
Then $R_1$:
$$\begin{array}{cccccccccc}
 &   &   &      &      &      &   &   & &\qquad[3]\\
 &   &   &      &b_{24}&      &   &   & &\qquad[2]\\
 &   &   &b_{23}&      &a_{45}&   &   & &\qquad[1]\\
 &   &y_2&      &   x_3&      &y_5&   & &\qquad[1]\\
 &x_1&   &   y_3&      &   y_4&   &y_5& &\qquad[1]\\
1&   &  2&      &     3&      &  4&   &5&\qquad[0]
\end{array}.$$
Next $R_3$:
$$\begin{array}{cccccccccc}
 &   &      &      &      &      &      &   & &\qquad[3]\\
 &   &      &      &b_{24}&      &      &   & &\qquad[2]\\
 &   &      &   y_3&      &   y_5&      &   & &\qquad[2]\\
 &   &b_{12}&      &a_{34}&      &a_{45}&   & &\qquad[1]\\
 &x_1&      &   y_3&      &   y_4&      &y_5& &\qquad[1]\\
1&   &     2&      &     3&      &     4&   &5&\qquad[0]
\end{array}.$$
Then $R_2$:
$$\begin{array}{cccccccccc}
 &   &      &      &      &      &      &   & &\qquad[3]\\
 &   &      &      &   y_4&      &      &   & &\qquad[3]\\
 &   &      &b_{13}&      &b_{35}&      &   & &\qquad[2]\\
 &   &b_{12}&      &a_{34}&      &a_{45}&   & &\qquad[1]\\
 &x_1&      &   y_3&      &   y_4&      &y_5& &\qquad[1]\\
1&   &     2&      &     3&      &     4&   &5&\qquad[0]
\end{array}.$$
Now $R_1$:
$$\begin{array}{cccccccccc}
 &   &      &      &      &      &      &   & &\qquad[4]\\
 &   &      &      &b_{14}&      &      &   & &\qquad[3]\\
 &   &      &b_{13}&      &b_{35}&      &   & &\qquad[2]\\
 &   &b_{12}&      &a_{34}&      &a_{45}&   & &\qquad[1]\\
 &x_1&      &   y_3&      &   y_4&      &y_5& &\qquad[1]\\
1&   &     2&      &     3&      &     4&   &5&\qquad[0]
\end{array}.$$
Next $R_4$:
$$\begin{array}{cccccccccc}
 &      &      &      &      &      &      &      & &\qquad[4]\\
 &      &      &      &b_{14}&      &      &      & &\qquad[3]\\
 &      &      &b_{13}&      &b_{35}&      &      & &\qquad[2]\\
 &      &   x_1&      &   y_4&      &   y_5&      & &\qquad[2]\\
 &b_{12}&      &b_{23}&      &a_{34}&      &a_{45}& &\qquad[1]\\
1&      &     2&      &     3&      &     4&      &5&\qquad[0]
\end{array}.$$
Then $R_3$:
$$\begin{array}{cccccccccc}
 &      &      &      &      &      &      &      & &\qquad[4]\\
 &      &      &      &b_{14}&      &      &      & &\qquad[3]\\
 &      &      &   x_1&      &   y_5&      &      & &\qquad[3]\\
 &      &b_{13}&      &b_{24}&      &b_{35}&      & &\qquad[2]\\
 &b_{12}&      &b_{23}&      &a_{34}&      &a_{45}& &\qquad[1]\\
1&      &     2&      &     3&      &     4&      &5&\qquad[0]
\end{array}.$$
Next $R_2$:
$$\begin{array}{cccccccccc}
 &      &      &      &      &      &      &      & &\qquad[4]\\
 &      &      &      &   x_1&      &      &      & &\qquad[3]\\
 &      &      &b_{14}&      &b_{25}&      &      & &\qquad[3]\\
 &      &b_{13}&      &b_{24}&      &b_{35}&      & &\qquad[2]\\
 &b_{12}&      &b_{23}&      &a_{34}&      &a_{45}& &\qquad[1]\\
1&      &     2&      &     3&      &     4&      &5&\qquad[0]
\end{array}.$$
Finally $R_1$:
$$\begin{array}{cccccccccc}
 &      &      &      &      &      &      &      & &\qquad[5]\\
 &      &      &      &a_{15}&      &      &      & &\qquad[4]\\
 &      &      &b_{14}&      &b_{25}&      &      & &\qquad[3]\\
 &      &b_{13}&      &b_{24}&      &b_{35}&      & &\qquad[2]\\
 &b_{12}&      &b_{23}&      &a_{34}&      &a_{45}& &\qquad[1]\\
1&      &     2&      &     3&      &     4&      &5&\qquad[0]
\end{array}.$$
Hence $\Phi$ takes the oriented complete monotone triangle (\ref{exOCMT})
to the tournament~(\ref{exTour}).

\section{Remarks}

We have defined $\Phi$ and $\Psi$ by using the operators $R_d$
and $L_d$ in a fixed order. However it is not hard to prove that
performing any achievable sequence of ${n\choose 2}$ raising
(or lowering) operations produces the same map.

\begin{prop}
Let $N=\frac12n(n-1)$ and suppose that $R_{u_j}$ $(1\le j\le N)$
are raising operators such that
$$\Phi'=R_{u_N}\circ R_{u_{N-1}}\circ\cdots\circ R_{u_2}\circ R_{u_1}$$
is a well-defined map from $\co_n$ to~$\ct_n$. Then
$\Phi'=\Phi$. Similarly if $L_{u_j}$ $(1\le j\le N)$
are lowering operators such that
$$\Psi'=L_{u_N}\circ L_{u_{N-1}}\circ\cdots\circ L_{u_2}\circ L_{u_1}$$
is a well-defined map from $\ct_n$ to~$\co_n$. Then
$\Psi'=\Psi$.
\end{prop}
\begin{pf}
We need only prove the assertion concerning raising operators, as that for
lowering operators follows by taking inverses. It is easy to see that
if $i\ne j$, $T$ is a triangle and both $R_i(T)$ and $R_j(T)$ exist,
then $|i-j|\ge2$ and $R_j\circ R_i(T)$ and $R_i\circ R_j(T)$ exist
and are equal. Also if $|i-j|\ge2$ and $R_i\circ R_j(T)$ exists
then $R_i(T)$ exists. Thus if $|i-j|\ge2$ we have
$R_i\circ R_j=R_j\circ R_i$ in the sense that if one side exists
so does the other and they are equal.

The definition of $\Phi$ is framed so that
$$\Phi=R_{v_N}\circ R_{v_{N-1}}\cdots\circ R_{v_2}\circ R_{v_1}$$
where for each~$j$, $v_j$ is the least$v$ such that
$$R_{v}\circ R_{v_{j-1}}\cdots\circ R_{v_2}\circ R_{v_1}$$
is well-defined on $\co_n$.

Suppose that $u_j\ne v_j$ for some~$j$, and take $j$ to be the
least such for which this holds. Then $u_j>v_j$.
Since $R_j$ increases the rank of row $j-1$ by 1, the sequences
of $u_i$s and of $v_i$s are permutations of each other. Hence for
some $k>j$ we have $u_k=v_j$ and take $k$ to be the least such. Let
$T'=R_{u_{j-1}}\circ\cdots \circ R_1(T)$ where $T$ is any element of~$\co_n$.
Then row $u_j$ of $T'$ is an $\cx$-row and all previous rows are $\cv$-rows.
When we apply a sequence of raising operators to $T'$ we cannot
apply an operator $R_w$ for $w<v_j$ until we have applied the operator
$T_{v_j}$ since the rows 1 to $v_j-1$ will remain as $\cv$-rows
and row $v_j$ as a $\cx$-row until this occurs. Also $T_{v_j+1}$
cannot be applied until $T_{v_j}$ has been since row $T_{v_j}$
will remain an $\cx$-row during this process.
Thus if $j\le i<k$ then $u_i\ge v_j+2$. Hence
\begin{eqnarray*}
R_{u_k}\circ R_{u_{k-1}}\circ\cdots\circ R_{u_{j+1}}\circ R_{u_j}
&=&R_{v_j}\circ R_{u_{k-1}}\circ\cdots\circ R_{u_{j+1}}\circ R_{u_j}\\
&=&R_{u_{k-1}}\circ\cdots\circ R_{u_{j+1}}\circ R_{u_j}\circ R_{v_j}.
\end{eqnarray*}
By using this substitution in the definition of~$\Phi'$, we
increase the value of~$j$. Repeating this process eventually gives
$\Phi'=\Phi$.
\end{pf}

We can define a notion of oriented alternating sign matrix by assigning
each $-1$ in a given alternating sign matrix an orientation (left or right).
This is an equivalent notation to oriented complete monotone triangle.
There is a natural correspondence between oriented alternating sign matrices
of order $n$ and domino tilings of the Aztec diamond of order $n-1$
(see~\cite{EKLP}). We have thus given an implicit bijection between
tournaments on vertex set $\{1,2,\ldots,n\}$ and
domino tilings of the Aztec diamond of order~$n-1$. In particular we
recover the result that there are $2^{n\choose 2}$ such tilings.

We can apply the argument in the proof of Theorem~\ref{main} to more general
triangles. Instead
of having the last row $1,2,\ldots,n$ choose any set $S$ of $n$ positive
integers and consider the oriented monotone triangles
$\co_S$ with bottom row consisting of the elements of $S$ in ascending
order. We wish to describe the set $\Phi(\co_S)$. Let $\ct_S$ be the set of
admissible triangles with bottom row consisting of the elements of~$S$,
and with each entry above the bottom row being an $a_{ij}$ or a $b_{ij}$.
The proof of Theorem~\ref{main} immediately adapts to show that
$\Phi$ and $\Psi$ are inverse bijections between $\co_S$ and~$\ct_S$.
It is immediate that $|\ct_S|=2^{n\choose 2}|\ct'_S|$ where $\ct'_S$
is the subset of $\ct_S$ consisting of triangles lacking entries~$b_{ij}$.
This is because each element of $\ct_S$ is derived from a unique
element of $\ct'_S$ by replacing a subset of its $a_{ij}$ entries by the
corresponding $b_{ij}$. Each triangle in $\ct'_S$ is determined by
the left values of the entries (recall that $a_{ij}$ has left value $i$).
By replacing each $a_{ij}$ in such a triangle by $i$ we get a bijection
between $\ct'_S$ and the set of \emph{strict monotone triangles}
with bottom row~$S$, that is the set of numerical triangles
in which each arrangement
$$\begin{array}{ccc}
 &j& \\
i& &k
\end{array}$$
satisfies $i\le j<k$. For instance there are 6 strict monotone triangles
with bottom row $S=\{1,2,4,5\}$ namely
$$
\begin{array}{ccccccc}
 & & &1& & & \\
 & &1& &2& & \\
 &1& &2& &4& \\
1& &2& &4& &5\\
\end{array},\ \
\begin{array}{ccccccc}
 & & &1& & & \\
 & &1& &3& & \\
 &1& &2& &4& \\
1& &2& &4& &5\\
\end{array},\ \
\begin{array}{ccccccc}
 & & &2& & & \\
 & &1& &3& & \\
 &1& &2& &4& \\
1& &2& &4& &5\\
\end{array},
$$
$$
\begin{array}{ccccccc}
 & & &1& & & \\
 & &1& &3& & \\
 &1& &3& &4& \\
1& &2& &4& &5\\
\end{array},\ \
\begin{array}{ccccccc}
 & & &2& & & \\
 & &1& &3& & \\
 &1& &3& &4& \\
1& &2& &4& &5\\
\end{array},\ \
\begin{array}{ccccccc}
 & & &2& & & \\
 & &2& &3& & \\
 &1& &3& &4& \\
1& &2& &4& &5\\
\end{array}.
$$
By adapting the proof in \cite[\S4]{EKLP}
one can show that if $S=\{a_1,a_2,\ldots,a_n\}$ with $a_1<a_2<\cdots<a_n$
then
$$|\ct'_S|=\prod_{1\le i<j\le n}\frac{a_j-a_i}{j-i}$$
and so
$$|\ct_S|=2^{{n\choose 2}}\prod_{1\le i<j\le n}\frac{a_j-a_i}{j-i}.$$

\end{document}